\documentclass[10pt]{article}
\usepackage{latexsym,amsfonts,amssymb,amsbsy}
\usepackage{graphicx}
\usepackage{enumerate}

  \newtheorem {lemma}      {Lemma}  [section]

  \newtheorem {theorem}    [lemma]  {Theorem}
  \newtheorem {remark}     [lemma]  {Remark}
  \newtheorem {corollary}  [lemma]  {Corollary}
  \newtheorem {proposition}[lemma]  {Proposition}
  
  \newtheorem {conjecture} [lemma]  {Conjecture}

\newcommand {\proof}[1][{}] {\hbox{\bf Proof #1.\ }}
\newcommand {\qed} {\null \hfill \rule{2mm}{2mm}}

\def\lover#1#2{\lower #1 ex\hbox{$#2$}}

\newcommand{\dogsa}[1]{{\Psi^{#1}_{#1}}}                            
\newcommand{\dogs}[2]{{\Psi^{#1}_{#2}}}                             
\newcommand{\dogston}{{\dogsa{1\ldots n}}}                          
\newcommand{\dogstonwok}{{\dogsa{1\ldots\widehat{k}\ldots n}}}      
\newcommand{\dogstona}[1]{{\dogsa{1\ldots {#1}}}}                   
\newcommand{\dogstonwoka}[1]{{\dogsa{1\ldots\widehat{k}\ldots {#1}}}} %
\newcommand{\twodots}[1]{{{1}..{#1}}}                               

\def\DJ{{D\kern-.8em\raise.27ex\hbox{-$\!$-}\kern.3em}}

\def\Z{{\mathbb Z}}
\def\Q{{\mathbb Q}}
\def\R{{\mathbb R}}
\def\C{{\mathbb C}}
\def\P{{\mathbb P}}

\title{ATIYAH-SUTCLIFFE CONJECTURES FOR ALMOST COLLINEAR CONFIGURATIONS AND SOME NEW
CONJECTURES FOR SYMMETRIC FUNCTIONS}
\author{Dragutin Svrtan
\\
\small Department of Mathematics, University of Zagreb,\\[-0.8ex]
\small Bijeni\v{c}ka cesta 30, 10000 Zagreb, Croatia,\\[-0.8ex]
\small \texttt{dsvrtan@math.hr}\\
Igor Urbiha
\\
\small Department of Mathematics, Polytechnic of Zagreb, University of Zagreb,\\[-0.8ex]
\small Konavoska 2, 10000 Zagreb, Croatia,\\[-0.8ex]
\small \texttt{urbiha@vtszg.hr}}

\date{
\small AMS Subject Classifications: 74H05, 11B37, 26A18, 05A15,
11Y55, 11Y65}

\begin{document}

\maketitle

%
%
%
%
\begin{abstract}

In 2001 Sir M.\ F.\ Atiyah formulated a conjecture $(C1)$ and
later with P.\ Sutcliffe two stronger conjectures $(C2)$ and
$(C3)$. These conjectures, inspired by physics (spin-statistics
theorem of quantum mechanics), are geometrically defined for any
configuration of points in the Euclidean three space. The
conjecture $(C1)$ is proved for $n=3,4$ and for general $n$ only
for some special
configurations (M.\ F.\ Atiyah, M.\ Eastwood and P.\ Norbury, D. {\DJ}jokovi\'{c}).\\

\indent In this paper we shall explain some new conjectures for
symmetric functions which imply $(C2)$ and $(C3)$ for almost
collinear configurations. Computations up to $n=6$ are performed
with a help of {\tt Maple} and J.\ Stembridge's package {\tt SF}
for symmetric functions. For $n=4$ the conjectures $(C2)$ and
$(C3)$ we have also verified for some infinite families of
tetrahedra.\ This is a joint work with I.~Urbiha.\\

 \indent Finally we mention that by minimizing  a
geometrically defined energy, figuring in these conjectures, one
gets a connection to some  complicated physical theories, such as
Skyrmions and Fullerenes.
\end{abstract}

\newcommand{\be}{\begin{equation}}
\newcommand{\ee}{\end{equation}}
\newcommand{\pr}{\partial}
\newcommand{\ie}{{\it i.e.\e }}
\newcommand{\bphi}{\mbox{\boldmath $\phi$}}
\newcommand{\bx}{{\bf x}}
\newcommand{\CP}{\C\P}

\section{Introduction on Geometric Energies}

In this Section we describe some geometric energies,
introduced by Atiyah.
To construct first geometric energy consider $n$ distinct
ordered points, ${\bf x}_i\in\R^3$ for $i=1,...,n$.
 For each pair $i\ne j$ define the unit vector
\be
{\bf v}_{ij}=\frac{\bx_j-\bx_i}{|\bx_j-\bx_i|}
\label{unitv}
\ee
giving the direction of the line joining $\bx_i$ to $\bx_j.$
Now let $t_{ij}\in\CP^1$ be the point on the Riemann sphere associated with
the unit vector ${\bf v}_{ij}$, via the identification
$\CP^1\cong S^2,$ realized as stereographic projection.
Next, set $p_i$ to be the polynomial in $t$ with
roots $t_{ij}$ ($j\ne i$), that is
\be
p_i=\alpha_i\prod_{j\ne i}(t-t_{ij})
\label{defp}
\ee
where $\alpha_i$ is a certain normalization coefficient.
In this way we have constructed $n$ polynomials which all have degree
$n-1,$ and so we may write
\[
p_i=\sum_{j=1}^n d_{ij}t^{j-1}.
\]
Finally, let $d$ be the
$n\times n$ matrix with entries $d_{ij},$ and let $D$ be its
determinant
\be
D({\bf x}_1,...,{\bf x}_n)=\mbox{det}\ d.
\ee
This geometrical construction is relevant to the Berry-Robbins
problem, which is concerned with specifying how a spin
basis varies as $n$ point particles move in space, and supplies
a solution provided it can be shown that $D$ is always non-zero.
For $n=2,3,4$ it can be proved that $D\ne 0$ (Eastwood and Norbury) and
numerical computations suggest that $|D|\ge 1$ for all $n,$
with the minimal value $|D|=1$ being attained by $n$ collinear points.

The geometric energy is the $n$-point energy defined by
\be
E=-\log |D|,
\label{geom}
\ee
so minimal energy configurations maximize the modulus of the determinant.

This energy is geometrical in the sense that it only depends on
the directions of the lines joining the points, so it is
translation, rotation and scale invariant. Remarkably, the minimal
energy configurations, studied numerically for all $n\le 32,$ are
essentially the same as those for the Thomson problem.

\section{Eastwood--Norbury formulas for Atiyah determinants}

In this section we first recall Eastwood--Norbury formula for
Atiyah determinant for three or four points in Euclidean
three--space. In the case $n=3$ the Atiyah determinant reads as
\[
\det M_3=d_3(r_{12},r_{13},r_{22})+8r_{12}r_{13}r_{22}
\]
where \
\[
d_3(a,b,c)=(a+b-c)(b+c-a)(c+a-b)
\]
and $r_{ij}$ $(1\leq i<j\leq 3)$ is the distance
between the $i^{\mbox{th}}$ and $j^{\mbox{th}}$ point.

In the case $n=4$ the Atiyah determinant $\det M_4$ has real part
given by a polynomial (with $248$ terms) as follows:
\begin{equation}
\label{E:1.0}
\Re (\det
M_4)=64r_{12}r_{13}r_{23}r_{14}r_{24}r_{34}-4d_3(r_{12}r_{34},r_{13}r_{24},r_{14}r_{23})+A_4+288V^2
\end{equation}
where
\[
A_4=\sum_{l=1}^{4}\left(\sum_{\scriptstyle (l\neq) i=1}^{4}
r_{li}((r_{lj}-r_{lk})^2-r_{jk}^2)\right)d_3(r_{ij},r_{ik},r_{jk})
\]
(here $\{j,k\}=\{1,2,3,4\}\setminus\{l,i\}$) and $V$ denotes the
volume of the tetrahedron with vertices our four points:
\[
\begin{array}{l@{\ =\ }l}
144V^2
&
-(r_{12}^4r_{34}^2+r_{12}^2r_{34}^4+\
\lover{-0.2}{\mbox{{\tiny four terms}}})+
(r_{12}^2r_{23}^2r_{34}^2+\ \lover{-0.2}{\mbox{{\tiny eleven
terms}}})-\\[1ex]
&
(r_{12}^2r_{23}^2r_{13}^2+\ \lover{-0.2}{\mbox{{\tiny
three terms}}}).
\end{array}
\]
We now state two formulas which will be used later:
\begin{enumerate}
    \item Alternative form of $A_4$:
\begin{equation}
\label{E:A4-a}
\begin{array}{ll}
A_4=
&\displaystyle
 \sum_{ {\scriptstyle l=1 } }^4
\left(( d_3(r_{il},r_{jl},r_{kl})+8r_{il}r_{jl}r_{kl}
+r_{il}(r_{il}^2-r_{jk}^2) +\right.\\
&
\hphantom{ \sum_{ {\scriptstyle l=1 } }^4}
r_{jl}(r_{jl}^2-r_{ik}^2)+\left. r_{kl}(r_{kl}^2-r_{ij}^2)\right)
d_3(r_{ij},r_{ik},r_{jk}),
\end{array}
\end{equation}
where for each $l$ we write $\{1,2,3,4\}\setminus\{l\}=\{i<j<k\}$.
    \item The sum of the second and the fourth term of (\ref{E:1.0}) can be rewritten as
\begin{equation}
\label{E:e0}
\begin{array}{l}
288V^2-4d_3(r_{12}r_{34},r_{13}r_{24},r_{14}r_{23})=\\
\hspace{0.7cm}\begin{array}{l}
=(r_{12}-r_{34})^2(r_{13}^2r_{24}^2+r_{14}^2r_{23}^2-r_{12}^2r_{34}^2)+\
\lover{-0.2}{\mbox{{\tiny two such terms}}}\ +\\
\hphantom{=}+4r_{12}r_{13}r_{23}r_{14}r_{24}r_{34}-\\
\hphantom{=}-r_{12}^2r_{13}^2r_{23}^2-r_{12}^2r_{14}^2r_{24}^2
-r_{13}^2r_{14}^2r_{34}^2-r_{23}^2r_{24}^2r_{34}^2.
\end{array}
\end{array}
\end{equation}
\end{enumerate}
It is well known that the quantity is always nonpositive.

The imaginary part $\Im (\det (M_4))$ of Atiyah determinant can
be written as a product of $144V^2$ with a polynomial (with
integer coefficients) having 369 terms.

The original Atiyah conjecture in our cases is equivalent to nonvanishing
of the determinants $\det (M_3)$ and $\det (M_4)$.

A stronger conjecture of Atiyah and Sutcliffe
(\cite{AS},Conjecture 2) states in our cases that $|\det (M_3)|$ $\geq$
$8r_{12}r_{13}r_{23}$ and $|\det (M_4)|\geq
64r_{12}r_{13}r_{23}r_{14}r_{24}r_{34}$.

From the formula (\ref{E:1.0}) above, with the help of the simple inequality
$d_3(a,b,c)$ $\leq$ $abc$ (for $a,b,c\geq 0$), Eastwood and Norbury got
"almost" the proof of the stronger conjecture by exhibiting the inequality
\[
\Re(\det M_4)\geq 60r_{12}r_{13}r_{23}r_{14}r_{24}r_{34}.
\]
To remove the word "almost" seems to be not so easy (at present
not yet done even for planar configuration of four points).

A third conjecture (stronger than the second) of Atiyah and
Sutcliffe (\cite{AS}, Conjecture 3) can be expressed, in the four
point case, in terms of polynomials in the edge lengths as
\[
|\det M_4|^2\geq
\prod_{\{i<j<k\}\subset\{1,2,3,4\}}(d_3(r_{ij},r_{ik},r_{jk})+8r_{ij}r_{ik}r_{jk})
\]
where the product runs over the four faces of the tetrahedron,
and the expressions for the left hand side are given explicitly.\\
(cf.\ {\tt
ftp://ftp.maths.adelaide.edu.au/meastwood/maple/points})

In this paper we study some infinite families of tetrahedra and
confirm the strongest conjecture of Atiyah and Sutcliffe for
several such infinite families.

\subsection{Atiyah--Sutcliffe conjecture for (vertically) upright tetrahedra (or pyramids)}

By upright tetrahedron we mean a tetrahedron with vertices $1$,
$2$, $3$, $4$ such that all edges from the vertex $4$ have equal
lengths. In this case we write
\[
r_{23}=a, r_{13}=b, r_{12}=c, r_{14}=r_{24}=r_{34}=d.
\]
Now we study the second and the third term of Eastwood--Norbury formula for $\Re (\det M_4)$:
\[
\begin{array}{l}
-4d_3(r_{12}r_{34},r_{13}r_{24},r_{14}r_{23})+A_4=-4d^3d_3(a,b,c)\\[2mm]
\hspace{1cm}+[c(b^2+2bd)+b(c^2+2cd)+d((b+c)^2-a^2)]a^2(2d-a)\\[2mm]
\hspace{1cm}+[c(a^2+2ad)+a(c^2+2cd)+d((c+a)^2-b^2)]b^2(2d-b)\\[2mm]
\hspace{1cm}+[b(a^2+2ad)+a(b^2+2bd)+d((a+b)^2-c^2)]c^2(2d-c)\\[2mm]
\hspace{1cm}+[d(4d^2-a^2)+d(4d^2-b^2)+d(4d^2-c^2)]d_3(a,b,c)\\[2mm]
=dB_1+B_2+B_3,
\end{array}
\]
where
{\small
\[
\begin{array}{l@{\!}l}
B_1 =& \hphantom{+}(8d^2-a^2-b^2-c^2)d_3(a,b,c),\\[2mm]
B_2 =& \hphantom{+}[a^2b^2c+a^2bc^2](2d-a)+
[a^2b^2c+ab^2c^2](2d-b)+[a^2bc^2+ab^2c^2](2d-c),\\[2mm]
B_3 =& \hphantom{+}(4bc+(b+c)^2-a^2)a^2(2d-a)+
(4ac+(a+c)^2-b^2)b^2(2d-b)+\\
& +(4ab+(a+b)^2-c^2)c^2(2d-c)
\end{array}
\]
}
The term $B_1$ could still be negative, so we start to employ the
geometric constraint on the edge length $d$:
\begin{center}
\begin{tabular}{l@{\ $=$\ }l}
$d\geq R$ & the circumradius of the base triangle $123$\\[2mm]
& $\displaystyle \frac{abc}{\sqrt{(a+b+c)d_3(a,b,c)}}$ (by Heron formula)
\end{tabular}
\end{center}
Now we split the term $B_2$ into two parts as follows:
\[
B_2 = dB_2'+B_2''
\]
where
{\small
\[
B_2'=
(a^2b^2c+a^2bc^2)\left(2-\frac{a}{R}\right)+
(a^2b^2c+ab^2c^2)\left(2-\frac{b}{R}\right)+
(a^2bc^2+ab^2c^2)\left(2-\frac{c}{R}\right)
\]
}
and
{\small
\[
B_2''=
(a^2b^2c+a^2bc^2)\left(\frac{ad}{R}-a\right)+
(a^2b^2c+ab^2c^2)\left(\frac{bd}{R}-b\right)+
(a^2bc^2+ab^2c^2)\left(\frac{cd}{R}-c\right).
\]
}
\begin{lemma}
\label{L:1}
For any nonnegative real numbers $a$, $b$, $c$ the inequality
\begin{equation}
\label{E:e1}
\begin{array}{l}
  [4abc(2abc+(a+b+c)(ab+ac+bc))-(a^2+b^2+c^2)D_3(a,b,c)]^2\geq\\[2mm]
\hspace{1cm}(a+b+c)^2(a^2b+ab^2+a^2c+ac^2+b^2c+bc^2)^2D_3(a,b,c)
\end{array}
\end{equation}
where $D_3(a,b,c)=(a+b+c)(a+b-c)(a-b+c)(-a+b+c)$ $(=(a+b+c)d_3(a,b,c))$, holds true.
\end{lemma}
\proof[of Lemma \ref{L:1}]
By assuming $a\geq b\geq c\geq 0$ and putting $a=b+h$ and $b=c+k$, $(h\geq 0, k\geq 0)$,
we get (by MAPLE, for example) a polynomial of degree $12$ with all coefficients positive
(between $1$ and $254930$). This is a computer verification of the Lemma \ref{L:1}. \qed

Now, by using $d\geq R$, we estimate
\[
\begin{array}{ll}
  B_1+B_2' &= [8d^2-(a^2+b^2+c^2)]d_3(a,b,c)+B_2'\\[2mm]
       &\geq  [8R^2-(a^2+b^2+c^2)]d_3(a,b,c)+B_2'\\[2mm]
           &= \hphantom{+}\frac{1}{a+b+c}(8a^2b^2c^2-(a+b+c)(a^2+b^2+c^2)d_3(a,b,c))+\\[1mm]
           &\hphantom{=}+4(a^2b^2c+a^2bc^2+ab^2c^2)-\\[1mm]
           &\hphantom{=}-(a^2b+ab^2+a^2c+ac^2+b^2c+bc^2)\sqrt{(a+b+c)d_3(a,b,c)}
\end{array}
\]
But this last quantity is nonnegative. This can be seen by taking square root on both sides of the inequality in
the Lemma \ref{L:1}. (Note that the quantity inside the brackets in the left hand side of the inequality
(\ref{E:e1}) is positive, which can be checked similarly.)
So, $B_1+B_2'\geq 0$. By putting together all what we have proved so far
we conclude that the Atiyah--Sutcliffe conjecture
\[
\Re(\det M_4)\geq 64r_{12}r_{13}r_{23}r_{14}r_{24}r_{34}
\]
is valid for upright tetrahedra.

\subsection{Atiyah--Sutcliffe conjectures for edge--tangential tetrahedra}

By edge--tangential tetrahedron we shall mean any tetrahedron for which there exists a sphere touching
all its edges (i.e.\ its 1--skeleton has an inscribed sphere.) For each $i$ from $1$ to $4$ we denote by
$t_i$ the length of the segment (lying on the tangent line) with one endpoint the vertex and the other
the point of contact of the tangent line with a sphere. Clearly for the distances $r_{ij}$ from
$i^{\mbox{th}}$ to $j^{\mbox{th}}$ point we have
\[
r_{ij}=t_i+t_j,\ (1\leq i<j\leq 4).
\]
Now we shall compute all the ingredients appearing in the Eastwood--Norbury formula for $\Re(\det M_4)$
in terms of elementary symmetric functions of the (tangential) variables $t_1, t_2, t_3, t_4$
(recall $e_1=t_1+t_2+t_3+t_4$, $e_2=t_1t_2+t_1t_3+t_1t_4+t_2t_3+t_2t_4+t_3t_4$,
$e_3=t_1t_2t_3+t_1t_2t_4+t_1t_3t_4+t_2t_3t_4$, $e_4=t_1t_2t_3t_4$).
\[
\begin{array}{l}
  64r_{12}r_{13}r_{23}r_{14}r_{24}r_{34} =
  \displaystyle 64\hspace{-2ex}\prod_{1\leq i<j\leq 4}\hspace{-2ex}(t_i+t_j)=
  64s_{3,2,1}=\\[6mm]
   \hphantom{64r_{12}r_{13}r_{23}r_{14}r_{24}r_{34}} =64\left|
\begin{array}{ccc}
  e_3 & e_4 & 0 \\
  e_1 & e_2 & e_3 \\
  0 & 1 & e_1
\end{array}
  \right|=
  64e_3e_2e_1-64e_4e_1^2-64e_3^2
\end{array}
\]
Here we have used Jacobi--Trudi formula for the triangular Schur function $s_{3,2,1}$
(see \cite{Jacobi}, (3.5)). Furthermore we have
\[
\begin{array}{r@{\ =\ }l}
  -4d_3(r_{12}r_{34},r_{13}r_{24},r_{14}r_{23}) & 128e_4e_2-32e_4e_1^2-32e_3^2 \\[2mm]
  288V^2 & 128e_4e_2-32e_3^2
\end{array}
\]
In order to compute $A_4$ we first compute, for fixed $l$ the following quantities
\[
\begin{array}{r@{\ =\ }l}
  d_3(r_{ij},r_{ik},r_{jk})  & 8t_it_jt_k\\
  \displaystyle \sum_{(l\neq) i=1}^4r_{li}((r_{lj}+r_{lk})^2-r_{jk}^2)
  & 4(3t_l(t_1+t_2+t_3+t_4)+(t_it_j+t_it_k+t_jt_k))t_l.
\end{array}
\]
Thus we get:
\[
A_4=32(3e_1^2+4e_2)e_4=96e_4e_1^2+128e_4e_2.
\]
Now we adjust terms in $\Re(\det M_4)$, in order to get shorter expression, as follows
\[
\begin{array}{rl}
  \Re(\det M_4) =& (64r_{12}r_{13}r_{23}r_{14}r_{24}r_{34}-2\cdot 288V^2)+\\
                & +(-4d_3(r_{12}r_{34},r_{13}r_{24},r_{14}r_{23})-288V^2)+A_4+4\cdot 288V^2 \\[2mm]
   =& (64e_3e_2e_1-64e_4e_1^2-256e_4e_2)+(-32e_4e_1^2)+\\
   & +(96e_4e_1^2+128e_4e_2)+ 4\cdot 288V^2\\[2mm]
   =&  64e_3e_2e_1-128e_4e_2+1152V^2\\[2mm]
   =&  64e_2(e_3e_1-2e_4)+1152V^2\\[2mm]
   =&  64e_2(2e_4+m_{211})+1152V^2,
\end{array}
\]
where $m_{211}=t_1^2t_2t_3+\cdots$ denotes the monomial symmetric function associated to the partition $(2,1,1)$.

In order to verify the third conjecture of Atiyah and Sutcliffe
\[
|\det M_4|^2\geq \prod_{\{i<j<k\}\subset\{1,2,3,4\}}(d_3(r_{ij},r_{ik},r_{jk})+8r_{ij}r_{ik}r_{jk})
\]
we note first that
\[
\begin{array}{l@{\ =\ }l}
  d_3(r_{ij},r_{ik},r_{jk})+8r_{ij}r_{ik}r_{jk} & (8t_it_jt_k+8(t_i+t_j)(t_i+t_k)(t_j+t_k)) \\[2mm]
   & 8(t_i+t_j+t_k)(t_it_j+t_it_k+t_jt_k)
\end{array}
\]
and state the following:
\begin{lemma}
\label{L:2}
For any nonnegative real numbers $t1,t_2,t_3,t_4\geq 0$ the following inequality
\begin{equation}
\label{E:e2}
\begin{array}{l}
(t_1t_2+t_1t_3+t_1t_4+t_2t_3+t_2t_4+t_3t_4)^2(2t_1t_2t_3t_4+m_{211}(t_1,t_2,t_3,t_4))^2\geq\\[2mm]
\geq \prod_{\{i<j<k\}\subset\{1,2,3,4\}}(t_i+t_j+t_k)(t_it_j+t_it_k+t_jt_k)
\end{array}
\end{equation}
holds true.
\end{lemma}
\proof[of Lemma \ref{L:2}]
The difference between the left hand side and the right hand side of the above inequality
(\ref{E:e2}), written in terms of monomial symmetric function is equal to
\[
\begin{array}{ll}
LHS-RHS=&m_{6321}+3m_{6222}+m_{543}+2m_{5421}+7m_{5322}+5m_{5331}+\\[2mm]
&+3m_{444}+7m_{4431}+8m_{4422}+8m_{4332}+3m_{3333}\geq 0
\end{array}
\]
\qed
\begin{remark}
One may think that the inequality in Lemma \ref{L:2} can be obtained as a product of two
simpler inequalities. This is not the case, because the following inequalities hold true:
{\small
\[
\begin{array}{rcl}
(t_1t_2+t_1t_3+t_1t_4+t_2t_3+t_2t_4+t_3t_4)^2
&\leq&
 \displaystyle\prod_{\{i<j<k\}\subset\{1,2,3,4\}}(t_i+t_j+t_k) \\[6mm]
(2t_1t_2t_3t_4+m_{211}(t_1,t_2,t_3,t_4))^2
&\geq&
 \displaystyle \prod_{\{i<j<k\}\subset\{1,2,3,4\}}(t_it_j+t_it_k+t_jt_k)  \\
    &
\end{array}
\]
}
\end{remark}
Now we continue with verification of the third conjecture of Atiyah and Sutcliffe for
edge tangential tetrahedron:
\[
\begin{array}{rl}
  |\det M_4|^2\geq & (\Re(\det M_4))^2\geq [64e_2(2e_4+m_{211})]^2 \\[2mm]
  \geq & \displaystyle 8^4\hspace{-5ex}\prod_{\{i<j<k\}\subset\{1,2,3,4\}}
  \hspace{-4ex}(t_i+t_j+t_k)(t_it_j+t_it_k+t_jt_k)\ \ \ (\mbox{by Lemma \ref{L:2}}) \\[6mm]
  = & \displaystyle\hspace{-5ex}\prod_{\{i<j<k\}\subset\{1,2,3,4\}}
  \hspace{-4ex}(d_3(r_{ij},r_{ik},r_{jk})+8r_{ij}r_{ik}r_{jk})
\end{array}
\]
so the strongest Atiyah--Sutcliffe conjecture is verified for edge--tangential tetrahedra.

\subsection{Atiyah--Sutcliffe conjectures for isosceles tetrahedra}

By an isosceles tetrahedron we shall mean a tetrahedron in which each pair of opposite edges are equal
(hence all triangular faces are congruent). A tetrahedron is isosceles iff the sum of the face angles
at each polyhedron vertex is $180^\circ$, and iff its insphere and circumsphere are concentric.
We now compute the Atiyah determinant for isosceles tetrahedra.

In our notation we have $r_{23}=r_{14}=a$, $r_{13}=r_{24}=b$, $r_{12}=r_{34}=c$.
By formula (\ref{E:A4-a}) we get immediately that
\[
-4d_3(r_{12}r_{34},r_{13}r_{24},r_{14}r_{23})+288V^2=0.
\]
By using our alternative formula (\ref{E:A4-a}) for the term $A_4$ in the Eastwood--Norbury formula
we get immediately
\[
\begin{array}{r@{\ =\ }l}
A_4&\displaystyle \sum_{l=1}^4\left(d_3(r_{il},r_{jl},r_{kl})+8r_{il}r_{jl}r_{kl}\right)
d_3(r_{ij},r_{ik},r_{jk})\\[6mm]
   & (2d_3(a,b,c)+8abc)^2
\end{array}
\]
where for each $l$ we write $\{1,2,3,4\}\setminus\{l\}=\{i<j<k\}$.\\
The real part of the Atiyah determinant is given by
\[
\begin{array}{r@{\ =\ }l}
\Re(\det(M_4)) & 64a^2b^2c^2+4(d_3(a,b,c)+8abc)d_3(a,b,c)\\
    & (2d_3(a,b,c)+8abc).
\end{array}
\]
Now
\[
|\det(M_4)|^2\geq \Re(\det(M_4))^2=(2d_3(a,b,c)+8abc)^4\geq (d_3(a,b,c)+8abc)^4
\]
verifies the third Atiyah--Sutcliffe conjecture for isosceles tetrahedra.


\subsection{Atiyah determinant for triangles and quadrilaterals via trigonometry}

Denote the three points $x_1$, $x_2$, $x_3$ simply by symbols $1, 2, 3$ and let $X$, $Y$ and $Z$ denote
the angles of the of the triangle at vertices $1$, $2$ and $3$ respectively. Then we can express
the Atiyah determinant $\det M_3=d_3(r_{12},r_{13},r_{23})+8r_{12}r_{13}r_{23}$ as follows
\[
\det M_3=4r_{12}r_{13}r_{23}\left(\cos^2\frac{X}{2}+\cos^2\frac{Y}{2}+\cos^2\frac{Z}{2}\right).
\]
This follows, by using cosine law and sum to product formula for cosine, from the following identity
\[
\begin{array}{l@{\ =\ }l}
d_3(a,b,c)+8abc&(a+b-c)(a-b+c)(-a+b+c)+8abc\\
               &a((b+c)^2-a^2)+b((c+a)^2-b^2)+c((a+b)^2-c^2).
\end{array}
\]
Now we shall translate the Eastwood--Norbury formula for (planar quadrilaterals) into a trigonometric form.
Denote the four points $x_1$, $x_2$, $x_3$, $x_4$ simply by symbols $1, 2, 3, 4$ and denote by
\[
(X^{(1)}, Y^{(1)}, Z^{(1)}),\ \
(X^{(2)}, Y^{(2)}, Z^{(2)}),\ \
(X^{(3)}, Y^{(3)}, Z^{(3)}),\ \
(X^{(4)}, Y^{(4)}, Z^{(4)})
\]
 the angles of the triangles $234$, $341$, $412$, $123$ in this cyclic order
 (i.e.\ the angle of a triangle $412$ at vertex $2$ is $Z^{(3)}$ etc.).

Next we denote by $c_l$, ($1\leq l\leq 4$), the sums of cosines squared of half-angles of the $l$--th triangle i.e.:
\[
c_l:=\cos^2\frac{X^{(l)}}{2}+\cos^2\frac{Y^{(l)}}{2}+\cos^2\frac{Z^{(l)}}{2},\ \ l=1,2,3,4.
\]
Similarly, we denote by $\widehat{c}_l$, ($1\leq l\leq 4$),
the sum of cosines squared of half-angles at the $l$--th vertex of our quadrilateral thus
\[
\begin{array}{c}
\displaystyle \widehat{c}_1=\cos^2\frac{Z^{(2)}}{2}+\cos^2\frac{Y^{(3)}}{2}+\cos^2\frac{X^{(4)}}{2}\\[3mm]
\displaystyle \widehat{c}_2=\cos^2\frac{Z^{(3)}}{2}+\cos^2\frac{Y^{(4)}}{2}+\cos^2\frac{X^{(1)}}{2}\\[3mm]
\displaystyle \widehat{c}_3=\cos^2\frac{Z^{(4)}}{2}+\cos^2\frac{Y^{(1)}}{2}+\cos^2\frac{X^{(2)}}{2}\\[3mm]
\displaystyle \widehat{c}_4=\cos^2\frac{Z^{(1)}}{2}+\cos^2\frac{Y^{(2)}}{2}+\cos^2\frac{X^{(3)}}{2}
\end{array}
\]
 Then the term
$A_4$ in the Eastwood--Norbury formula can be rewritten as
\[
\begin{array}{r@{\ =\ }l}
  A_4 & \displaystyle\sum_{l=1}^4(4r_{li}r_{lj}r_{lk}\widehat{c}_l)\cdot 4r_{ij}r_{ik}r_{jk}(c_l-2) \\
    & \displaystyle16r_{12}r_{13}r_{23}r_{14}r_{24}r_{34}\sum_{l=1}^4\widehat{c}_l(c_l-2). \\
\end{array}
\]
where for each $l$ we write $\{1,2,3,4\}\setminus\{l\}=\{i<j<k\}$.\\

In order to rewrite the term $-4d_3(r_{12}r_{34},r_{13}r_{24},r_{14}r_{23})$ into a trigonometric form
we recall a theorem
of M\" obius (\cite{Moebius}) which claims that for any quadrilateral $1234$ in a plane the products
$r_{12}r_{34}$, $r_{13}r_{24}$ and $r_{14}r_{23}$ are proportional to the sides of a triangle
whose angles are the differences of angles in the quadrilateral $1234$:
\[
\begin{array}{c@{\ =\ }c@{\ =\ }l}
  X & \sphericalangle 243 - \sphericalangle 213 & X^{(4)}-Z^{(1)} \\
  Y & \sphericalangle 341 - \sphericalangle 321 & Y^{(2)}-(-Y^{(4)}) \\
  Z & \sphericalangle 142 - \sphericalangle 132 & -X^{(3)}+Z^{(4)}
\end{array}
\]

Thus
\[
-4d_3(r_{12}r_{34},r_{13}r_{24},r_{14}r_{23})=-16r_{12}r_{13}r_{23}r_{14}r_{24}r_{34}(c-2)
\]
where
\[
c=\cos^2\frac{X}{2}+\cos^2\frac{Y}{2}+\cos^2\frac{Z}{2}.
\]
Thus we have obtained a trigonometric formula for Atiyah determinant of quadrilaterals
\[
\begin{array}{l@{\ =\ }l}
\Re(\det M_4)&\displaystyle\prod_{1\leq i<j\leq 4}r_{ij}\left(64-16(c-2)+16\sum_{l=1}^4\widehat{c}_l(c_l-2)\right)\\[4mm]
&\displaystyle 16\prod_{1\leq i<j\leq 4}r_{ij}\left(6-c+\sum_{l=1}^4\widehat{c}_l(c_l-2)\right)
\end{array}
\]
Now we shall verify Atiyah--Sutcliffe conjecture for cyclic quadrilaterals. In this case, by
a well known Ptolemy's theorem, we see that
\[
-4d_3(r_{12}r_{34},r_{13}r_{24},r_{14}r_{23})=0\ \ (\Leftrightarrow c=2)
\]
By using the equality of angles $Z^{(2)}=X^{(1)}$, $Z^{(3)}=X^{(2)}$, $Z^{(4)}=X^{(3)}$, $Z^{(1)}=X^{(4)}$
and $Y^{(1)}+Y^{(3)}=\pi=Y^{(2)}+Y^{(4)}$ (angles with vertex on a circle's circumference with the same
endpoints are equal or suplement of each other)we obtain
\[
\begin{array}{c}
\displaystyle \widehat{c}_1=\cos^2\frac{X^{(1)}}{2}+\cos^2\frac{Y^{(1)}}{2}+\cos^2\frac{Z^{(1)}}{2}=c_1-\cos Y^{(1)},\\[3mm]
\displaystyle \widehat{c}_2=\cos^2\frac{X^{(2)}}{2}+\cos^2\frac{Y^{(2)}}{2}+\cos^2\frac{Z^{(2)}}{2}=c_2-\cos Y^{(2)},\\[3mm]
\displaystyle \widehat{c}_3=\cos^2\frac{X^{(3)}}{2}+\cos^2\frac{Y^{(3)}}{2}+\cos^2\frac{Z^{(3)}}{2}=c_3-\cos Y^{(3)},\\[3mm]
\displaystyle \widehat{c}_4=\cos^2\frac{X^{(4)}}{2}+\cos^2\frac{Y^{(4)}}{2}+\cos^2\frac{Z^{(4)}}{2}=c_4-\cos Y^{(4)}.
\end{array}
\]
Now we have
\[
\begin{array}{rl}
\Re (\det M_4)=&
    \displaystyle \left(\prod_{1\leq i<j\leq 4}r_{ij} \right)
    \left(64+16\sum_{l=1}^4(c_l-\cos Y^{(l)})(c_l-2)\right)\\[6mm]
\geq&
\displaystyle \left(\prod_{1\leq i<j\leq 4}r_{ij} \right)
    \left(64+16\sum_{l=1}^4(c_l-1)(c_l-2)\right)\\[6mm]
\end{array}
\]
(here we have used that $2\leq c_l (\leq\frac{9}{4})$ for each $l=1,2,3,4$)
\[
\begin{array}{rl}
\hphantom{\Re (\det M_4)}\geq&
\displaystyle \left(\prod_{1\leq i<j\leq 4}r_{ij} \right)
    \left(64+16\sum_{l=1}^4(c_l-2)+16\sum_{l=1}^4(c_l-2)^2\right)\\[6mm]
\geq&
\displaystyle \left(\prod_{1\leq i<j\leq 4}r_{ij} \right)
    \left(64+16\sum_{l=1}^4(c_l-2)+4\left(\sum_{l=1}^4(c_l-2)\right)^2\right)
\end{array}
\]
(by quadratic--arithmetic inequality)
\[
\begin{array}{rl}
\hphantom{\Re (\det M_4)}=&
\displaystyle \left(\prod_{1\leq i<j\leq 4}r_{ij} \right)
    \left(\left(8+\sum_{l=1}^4(c_l-2)\right)^2+3\left(\sum_{l=1}^4(c_l-2)\right)^2\right)\\[6mm]
=&
\displaystyle \left(\prod_{1\leq i<j\leq 4}r_{ij} \right)
    \left(\left(\sum_{l=1}^4c_l\right)^2+\left(3\sum_{l=1}^4(c_l-2)\right)^2\right)\\[6mm]
\geq&
\displaystyle \left(\prod_{1\leq i<j\leq 4}r_{ij} \right)
    \left(\sum_{l=1}^4c_l\right)^2\geq 16\sqrt{c_1c_2c_3c_4}\prod_{1\leq i<j\leq 4}r_{ij}
\end{array}
\]
by A--G inequality.
Finally,
\[
\begin{array}{rl}
|\det M_4|^2=&\displaystyle
|\Re(\det M_4)|^2\geq 4^4c_1c_2c_3c_4\prod_{1\leq i<j\leq 4}r_{ij}^2\\
=&\displaystyle
    \prod_{l=1}^4(4r_{ij}r_{ik}r_{jk}c_l)
=
    \prod_{l=1}^4(d_3(r_{ij},r_{ik},r_{jk})+8r_{ij}r_{ik}r_{jk})
\end{array}
\]
where for each $l$ we write $\{1,2,3,4\}\setminus\{l\}=\{i<j<k\}$.
This finishes verification of Atiyah--Sutcliffe conjectures for cyclic quadrilaterals.



\section{Almost collinear configurations. {\DJ}okovi\' c's approach}

\subsection{Type (A) configurations}

By a type (A) configurations of $N$ points $x_1,\ldots,x_N$ we shall
mean the case when $N-1$ of the points $x_1,\ldots,x_N$ are
collinear. Set $n=N-1$. In (\cite{DZ}) {\DJ}okovi\' c has proved,
for configurations of type (A), both the Atiyah conjecture
(Theorem 2.1) and the first Atiyah--Sutcliffe conjecture (Theorem
3.1). By using Cartesian coordinates, with $x_i=(a_i,0)$,
$a_1<a_2<\cdots <a_n$ and $x_N=x_{n+1}=(0,b)'$ (with $b=1$), the normalized
Atiyah matrix $M_{n+1}=M_{n+1}(\lambda_1,\ldots,\lambda_n)$
(denoted by $P$ in \cite{DZ} when $b=-1$) is given by
\[
M_{n+1}= \left[
\begin{array}{cccccc}
  1 & \lambda_1 & 0 & \cdots & 0 & 0 \\
  0 & 1 & \lambda_2 & \cdots & 0 & 0 \\
  0 & 0 & 1 & & 0 & 0\\
  \vdots & \vdots & \vdots & \ddots & \vdots & \vdots \\
  0 & 0 &  &  & 1 & \lambda_n \\
  (-1)^{n}e_{n} & (-1)^{n-1}e_{n-1} & \cdots & \cdots & -e_1 & 1
\end{array}
\right]
\]
where $\lambda_1=a_1+\sqrt{a_1^2+b^2}$ $<$
$\lambda_2=a_2+\sqrt{a_2^2+b^2}$ $<\cdots<$
$\lambda_n=a_n+\sqrt{a_n^2+b^2}$ (with $b=1$) are positive real numbers and
where $e_k=e_k(\lambda_1,\ldots,\lambda_n)$, $1\leq k\leq n$, is
the $k$--th elementary symmetric function of
$\lambda_1,\lambda_2,\ldots,\lambda_n$. Its determinant satisfies the inequality
\[
\begin{array}{rl}
  \det(M_4) =& 1+\lambda_ne_1+\lambda_n\lambda_{n-1}e_2+\cdots+\lambda_n\lambda_{n-1}\cdots\lambda_1e_n\\
  \geq & 1+e_1(\lambda_1^2,\ldots,\lambda_n^2)+e_2(\lambda_1^2,\ldots,\lambda_n^2)+
  \cdots+e_n(\lambda_1^2,\ldots,\lambda_n^2) \\
  = & \prod_{i=1}^n(1+\lambda_i^2) \\
\end{array}
\]
equivalent to the first Atiyah--Sutcliffe conjecture (\cite{AS},Conjecture 2). The second
Atiyah--Sutcliffe conjecture (\cite{AS},Conjecture 3) for configurations of type (A) is equivalent
to the following inequality
\begin{equation}
\label{E:2.1}
[\det M_{n+1}(\lambda_1,\ldots,\lambda_n)]^{n-1}\geq
\prod_{k=1}^n\det M_n(\lambda_1,\ldots,\lambda_{k-1},\lambda_{k+1},\ldots,\lambda_n)
\end{equation}
For $n=2$ this inequality takes the form
\[
1+\lambda_2e_1(\lambda_1,\lambda_2)+\lambda_1\lambda_2e_2(\lambda_1,\lambda_2)
\geq
(1+\lambda_2e_1(\lambda_2))(1+\lambda_1e_1(\lambda_1)
\]
i.e.
\begin{equation}
\label{E:2.2}
1+\lambda_2e_1(\lambda_1,\lambda_2)+\lambda_1\lambda_2e_2(\lambda_1,\lambda_2)
\geq
(1+\lambda_2^2)(1+\lambda_1^2).
\end{equation}
This reduces to $(\lambda_2-\lambda_1)\lambda_1\geq0$, so it is true.

Even for $n=3$ the inequality (\ref{E:2.1}) is quite messy thanks
to nonsymmetric character of both sides. Knowing that sometimes it
is easier to solve a more general problem we followed that path
(although we didn't solve the problem in full generality). So let
us start with the case $n=2$. If we look at the following
inequality
\[
1+X_1(\xi_1+\xi_2)+X_1X_2\xi_1\xi_2\geq (1+X_1\xi_1)(2+X_2\xi_2)
\]
which is clearly true if $X_1\geq X_2\geq 0$ and $\xi_1,\xi_2\geq 0$ we obtain the inequality
(\ref{E:2.2}) simply by a specialization $X_1=\xi_1=\lambda_2$, $X_2=\xi_2=\lambda_1$.
So we proceed as follows:


Let $\xi_1,\ldots,\xi_n, X_1,\ldots, X_n, n\geq 1$  be two sets of commuting indeterminates. For any $l, 1\leq
l\leq n$ and any sequences $1\leq i_1\leq\cdots\leq i_l\leq n, 1\leq j_1\leq\cdots\leq j_l\leq n$ we define
polynomials $\Psi^I_J=\Psi^{i_1,\ldots,i_l}_{j_1,\ldots,j_l}\in\Q[\xi_1,\ldots,\xi_n, X_1,\ldots, X_n]$ as
follows:
\[
\Psi^I_J:=\sum_{k=0}^le_k(\xi_{j_1},\xi_{j_2},\ldots,\xi_{j_l})X_{i_1}X_{i_2}\cdots
X_{i_k},\ (l\geq 1),\ \Psi^\emptyset_\emptyset:=1\ (j=0)
\]
where $e_k$ is the $k$-th elementary symmetric function.

In particular we have
\[
\begin{array}{ll}
\Psi^i_j & =1+\xi_jX_i,\\[2mm]
\Psi^{i_1i_2}_{j_1j_2}& = 1+(\xi_{j_1}+\xi_{j_2})X_{i_1}+\xi_{j_1}\xi_{j_2}X_{i_1}X_{i_2},\\[2mm]
\Psi^{i_1i_2i_3}_{j_1j_2j_3} & =1+(\xi_{j_1}+\xi_{j_2}+\xi_{j_3})X_{i_1}+
(\xi_{j_1}\xi_{j_2}+\xi_{j_1}\xi_{j_3}+\xi_{j_2}\xi_{j_3})X_{i_1}X_{i_2}+\\
& \hphantom{= 1}+\xi_{j_1}\xi_{j_2}\xi_{j_3}X_{i_1}X_{i_2}X_{i_3},\\
\mbox{ etc.}
\end{array}
\]

The polynomials $\Psi^I_J$ are symmetric w.r.t.\
$\xi_{j_1},\xi_{j_2},\ldots,\xi_{j_l}$, but nonsymmetric w.r.t.\
$X_{i_1},X_{i_2},\ldots,X_{i_l}$. These polynomials when
restricted to nonnegative arguments such that $X_{i_1}\geq
X_{i_2}\geq \ldots\geq X_{i_l}\geq 0$,
$\xi_{j_1},\xi_{j_2},\ldots,\xi_{j_l}\geq 0$ obey some intriguing
inequalities which are not yet proved in full generality .In turn
they generalize some special cases of not yet proven conjectures
of Atiyah and Sutcliffe on configurations of points in three
dimensional Euclidean space (the former inequalities are not yet
proven even for the case of $n=4$ points).

Let us now formulate a conjecture which implies the strongest
Atiyah--Sutcliffe's conjecture for almost collinear configurations
of points (all but one point are collinear, called type(A) in
\cite{DZ}).

Our conjecture reads as follows:

\begin{conjecture}
\label{C:1} For any $n\geq 1$, let $X_{i_1}\geq X_{i_2}\geq
\ldots\geq X_{i_l}$ $\geq$ $0$,
$\xi_{j_1},\xi_{j_2},\ldots,\xi_{j_l}\geq 0$, be any nonnegative
real numbers. Then
\[
\left(
\Psi^{12\cdots n}_{12\cdots n}
\right)^{n-1}
\geq
\prod_{k=1}^{n}\Psi^{12\cdots \hat{k}\cdots n}_{12\cdots \hat{k}\cdots n}
\]
where $12\cdots \hat{k}\cdots n$ denotes the sequence $12\cdots (k-1)(k+1)\cdots n$. The equality obviously holds
true iff $X_1=X_2=\cdots =X_n$.
\end{conjecture}

To illustrate the Conjecture (\ref{C:1}) we consider first the cases $n=2$ and $n=3$.

\begin{description}
    \item[Case $n=2$:] We have
    \[
        \begin{array}{l@{}l}
          \Psi^{12}_{12} & = 1+(\xi_1+\xi_2)X_1+\xi_1\xi_2X_1X_2=\\[2mm]
                         & = 1+\xi_1X_1+\xi_2X_2+\xi_1\xi_2X_1X_2+(X_1-X_2)\xi_2=\\[2mm]
                         & = (1+\xi_1X_1)(1+\xi_2X_2)+\xi_2(X_1-X_2)\geq\\[2mm]
                         & \geq (1+\xi_1X_1)(1+\xi_2X_2)=\Psi^1_1\Psi^2_2.
        \end{array}
    \]
    \item[Case $n=3$:] We first write $\Psi^{123}_{123}$ in two different ways:
    \[
        \Psi^{123}_{123}=\xi_2(X_1-X_2)+\widehat{\Psi}^{123}_{1\underline{2}3} \mbox{\ \ \  and\ \ \  }
        \Psi^{123}_{123}=\xi_3(X_1-X_2)+\widehat{\Psi}^{123}_{12\underline{3}}.
    \]
    Note that $\widehat{\Psi}^{123}_{1\underline{2}3}$ is obtained from $\Psi^{123}_{123}$ by replacing
    the linear term $\xi_2X_1$ by $\xi_2X_2$, hence all its coefficients are nonnegative.

    The left hand side of the Conjecture (\ref{C:1}) $L_3$ can be rewritten as follows:
    \[
\begin{array}{l@{}l}
  L_3=(\Psi^{123}_{123})^2 & = (\xi_2(X_1-X_2)+\widehat{\Psi}^{123}_{1\underline{2}3})\Psi^{123}_{123}\\[2mm]
                   & = \xi_2(X_1-X_2)\Psi^{123}_{123}+\widehat{\Psi}^{123}_{1\underline{2}3}\Psi^{123}_{123}\\[2mm]
    & = \xi_2(X_1-X_2)\Psi^{123}_{123}+
        \widehat{\Psi}^{123}_{1\underline{2}3}(\xi_3(X_1-X_2)+\widehat{\Psi}^{123}_{12\underline{3}})\\[2mm]
    & = L'_3(X_1-X_2)+\widehat{\Psi}^{123}_{1\underline{2}3}\widehat{\Psi}^{123}_{12\underline{3}}
\end{array}
    \]
    where $L'_3=\xi_2\Psi^{123}_{123}+\xi_3\widehat{\Psi}^{123}_{1\underline{2}3}$ is a positive polynomial.

    Now we have
    \[
    L_3\geq \widehat{L}_3:=\widehat{\Psi}^{123}_{1\underline{2}3}\widehat{\Psi}^{123}_{12\underline{3}}.
    \]
    By using the formula
    \[
    \widehat{\Psi}^{123}_{1\underline{2}3}=\Psi^{12}_{13}+\xi_2X_2\Psi^{13}_{13}=
        (\Psi^2_{2}-1)\Psi^{13}_{13}+\Psi^{12}_{13}
    \]
    we can rewrite $\widehat{L}_3$ as
    \[
\begin{array}{l@{}l}
  \widehat{L}_3 & =\left[ (\Psi^{12}_{13}-\Psi^{13}_{13})+
                           \Psi^2_2\Psi^{13}_{13}\right]\widehat{\Psi}^{123}_{12\underline{3}}\\[2mm]
   & = \xi_1\xi_3X_1(X_2-X_3)\widehat{\Psi}^{123}_{12\underline{3}}+
       \Psi^{13}_{13}(\Psi^2_2\widehat{\Psi}^{123}_{12\underline{3}})
\end{array}
    \]
    The last term in parenthesis can be written as
    \[
\begin{array}{l@{}l}
  \Psi^2_2\widehat{\Psi}^{123}_{12\underline{3}} & = \Psi^{12}_{12}\Psi^{23}_{23}+
                                                     \Psi^1_2(\Psi^{22}_{23}-\Psi^{23}_{23})\\[2mm]
    & = \Psi^{12}_{12}\Psi^{23}_{23}+\xi_2\xi_3X_2(X_2-X_3)\Psi^1_2,
\end{array}
    \]
    so we get
    \[
    \widehat{L}_3=L''_3(X_2-X_3)+\Psi^{12}_{12}\Psi^{13}_{13}\Psi^{23}_{23}
    \]
    where $L''_3$ denotes the positive polynomial
    \[
    L''_3=\xi_1\xi_3X_1\widehat{\Psi}^{123}_{12\underline{3}}+\xi_2\xi_3X_2\Psi^1_2\Psi^{13}_{13}.
    \]
    We now have an explicit formula for $L_3$:
    \[
    L_3=L'_3(X_1-X_2)+L''_3(X_2-X_3)+\Psi^{12}_{12}\Psi^{13}_{13}\Psi^{23}_{23}
    \]
    with $L'_3, L''_3$ positive polynomials, which together with $X_1\geq X_2\geq X_3 (\geq 0)$ implies that
    \[
    L_3\geq R_3:=\Psi^{12}_{12}\Psi^{13}_{13}\Psi^{23}_{23}
    \]
    and the Conjecture (\ref{C:1}) ($n=3$) is proved.
\end{description}

In fact we have proven an instance $n=3$ $\widehat{L}_3\geq R_3$
of a stronger conjecture which we are going to formulate now. Let
$2\leq k\leq n$. We define the modified polynomials
$\widehat{\Psi}^{12\ldots k\ldots n}_{12\ldots \underline{k}\ldots
n}$ as follows:
\[
    \widehat{\Psi}^{12\ldots k\ldots n}_{12\ldots \underline{k}\ldots n}:=
    \xi_k(X_2-X_1)+\Psi^{12\ldots n}_{12\ldots n}
\]
obtained from $\Psi^{12\ldots n}_{12\ldots n}$ by replacing only
one term $\xi_kX_1$ by $\xi_kX_2$, hence $\widehat{\Psi}^{12\ldots
k\ldots n}_{12\ldots \underline{k}\ldots n}$ are still positive.
Let us introduce the following notation:
\[
    \widehat{L}_n:=\prod_{k=2}^{n}\widehat{\Psi}^{12\ldots k\ldots n}_{12\ldots \underline{k}\ldots n}\ ;\ \
    R_n:=\prod_{k=1}^n\Psi^{12\ldots \hat{k}\ldots n}_{12\ldots \hat{k}\ldots n}.
\]
Then clearly $L_n:=(\Psi^{12\ldots n}_{12\ldots n})^{n-1}\geq \widehat{L}_n$.
Now our stronger conjecture reads as
\begin{conjecture}
\label{C:2}
\[\displaystyle \widehat{L}_n\geq R_n\  (n\geq 1)\]
with equality iff $X_2=X_3=\cdots =X_n$.
\end{conjecture}
More generally, we conjecture that the difference $\widehat{L}_n-R_n$ is a polynomial in the
differences $X_2-X_3$, $X_3-X_4$, $\ldots$, $X_{n-1}-X_{n}$ with coefficients in
$\Z_{\geq 0}[X_1,\ldots,X_n,\xi_1,\ldots,\xi_n]$.
\begin{proposition}
\label{P:1}
\[L_n=L'_n(X_1-X_2)+\widehat{L}_n\]
for some positive polynomial $L'_n$.
\end{proposition}
\proof[of Proposition \ref{P:1}]
\[
\begin{array}{l}
 L_n=(\Psi^{12\cdots n}_{12\cdots n})^{n-1}=(\xi_2(X_1-X_2)+\widehat{\Psi}^{12\cdots n}_{1\underline{2}\cdots n})
        (\Psi^{12\cdots n}_{12\cdots n})^{n-2}\\[2mm]
=\xi_2(X_1-X_2)(\Psi^{12\cdots n}_{12\cdots n})^{n-2}+
        \widehat{\Psi}^{12\cdots n}_{1\underline{2}\cdots n}(\xi_3(X_1-X_2)+
        \widehat{\Psi}^{123\cdots n}_{12\underline{3}\cdots n})(\Psi^{12\cdots n}_{12\cdots n})^{n-3}\\[2mm]
=\xi_2(X_1-X_2)(\Psi^{12\cdots n}_{12\cdots n})^{n-2}+
\xi_3(X_1-X_2)\widehat{\Psi}^{12\cdots n}_{1\underline{2}\cdots n}(\Psi^{12\cdots n}_{12\cdots n})^{n-3}+\\
\hphantom{=}+\widehat{\Psi}^{12\cdots n}_{1\underline{2}\cdots n}\widehat{\Psi}^{123\cdots n}_{12\underline{3}\cdots n}
(\Psi^{12\cdots n}_{12\cdots n})^{n-3}\\
\vdots     \\
=(\sum_{k=1}^{n-1}\xi_{k+1}(\prod_{j=2}^k\widehat{\Psi}^{12\ldots j\ldots n}_{12\ldots \underline{j}\ldots n})
(\Psi^{12\ldots n}_{12\ldots n})^{n-k})(X_1-X_2)+
\prod_{j=2}^{n}\widehat{\Psi}^{12\ldots j\ldots n}_{12\ldots \underline{j}\ldots n}.
\end{array}
\]
\qed

Now we turn to study the quotient
\[
\frac{L_n}{R_n}=\frac{(\Psi^{1\ldots n}_{1\ldots n})^{n-1}}
{\displaystyle \prod_{k=1}^{n}\Psi^{1\ldots \widehat{k}\ldots n}_{1\ldots \widehat{k}\ldots n}}
\]
by studying the growth behaviour of quotients of its factors
$\Psi^{1\ldots n}_{1\ldots n}/
{\Psi^{1\ldots \widehat{k}\ldots n}_{1\ldots \widehat{k}\ldots n}}$
w.r.t.\ any of its arguments $X_r$, $1\leq r\leq n$.

\newpage
In the following theorem we obtain an explicit formula for the
numerators of the logarithmic derivatives w.r.t.\ $X_r, (1\leq r\leq n,
r\neq k)$ of the quantities $\displaystyle {\dogston/\dogstonwok}$.
From this formulas we get some monotonicity properties which enable us to state
some new (refined) conjectures later on.
\begin{theorem}
\label{T:1}
Let
  \begin{equation}
  \label{Eq:0r}
  \Delta_r:=\partial_{X_r}\dogs{1\ldots n}{1\ldots n}\cdot\dogsa{1\ldots\widehat{k}\ldots n}-
            \dogsa{1\ldots n}\cdot\partial_{X_r}\dogsa{1\ldots\widehat{k}\ldots n},
            \ \ (1\leq r\leq n).
  \end{equation}

Then we have the following explicit formulas
\begin{description}
\item[(i)] for any $r$, $1\leq r<k(\leq n)$ we have
\[
    \begin{array}{ll}
    \Delta_r= &\displaystyle \xi_k\!\!\!\!\!\!\!\!\!\!\sum_{0\leq i< r\leq j\leq n}\!\!\!\!\!\!
    s^{(k)}_{(2^i1^{j-i-1})}X_1^2\cdots X_i^2X_{i+1}\cdots \widehat{X}_k\cdots X_j+\\
    &\displaystyle+\!\!\!\!\!\!\!\!\!\!\sum_{0\leq i<r,k\leq j<n}\!\!\!\!\!\!\!\!
    e_ie_j^{(k)}X_1^2\cdots X_i^2X_{i+1}\cdots \widehat{X}_r\cdots \widehat{X}_k\cdots X_j(X_k-X_{j+1})
    \end{array}
  \]
\item[(ii)] for any $r$, $(1\leq )k<r\leq n$ we have
\[
\begin{array}{ll}
\Delta_r=&\displaystyle -\left(
\sum_{0\leq i< r\leq j\leq n}\!\!\!\!\!\!
    s^{(k)}_{(2^i1^{j-i-1})}X_1^2\cdots X_i^2X_{i+1}\cdots \widehat{X}_k\cdots \widehat{X}_r\cdots X_j+\right.\\
    &\displaystyle\left. +\!\!\!\!\!\!\!\!\!\!\sum_{0\leq i<k,r\leq j<n}\!\!\!\!\!\!\!\!
    e_i^{(k)}e_jX_1^2\cdots X_i^2X_{i+1}\cdots \widehat{X}_k\cdots \widehat{X}_r\cdots X_j(X_{j+1}-X_k)
\right)
    \end{array}
\]
where $s^{(k)}_\lambda$ denotes the $\lambda$--th Schur function of
$\xi_1,\ldots ,\xi_{k-1},\xi_{k+1},\ldots,\xi_n$ ($\xi_k$ omitted).
\end{description}

\end{theorem}
\proof[of Theorem \ref{T:1}]

  {\bf (i)} For any $r$, $1\leq r<k(\leq n)$ we find explicitly a formula for
  as follows. We shall use notations $X_\twodots{i}:=X_1X_2\cdots X_i$,
  for multilinear monomials and $e_i:=e_i(\xi_1,\ldots,\xi_n)$,
  $e_i^{(k)}=e_i(\xi_1,\ldots,\widehat{\xi_k},\ldots\xi_n)$ for the elementary
  symmetric functions (here $k$ is fixed). Then we can rewrite our basic quantities
  \begin{equation}
  \label{Eq:1r}
  \dogsa{1\ldots n}:=\sum_{i=0}^{n}e_iX_\twodots{i}
  \end{equation}
  \begin{equation}
  \label{Eq:2r}
  \begin{array}{lrl}
  \dogsa{1\ldots\widehat{k}\ldots n}&:=&\displaystyle\sum_{i=0}^{k-1}e_i^{(k)}X_\twodots{i}+
  \frac{1}{X_k}\sum_{i=0}^{n-1}e_i^{(k)}X_\twodots{i+1}=\\[4mm]
    &=&\displaystyle\sum_{i=0}^{n-1}e_i^{(k)}X_\twodots{i}+
    \frac{1}{X_k}\sum_{i=0}^{n-1}e_i^{(k)}X_\twodots{i}(X_{i+1}-X_k)
  \end{array}
  \end{equation}
  For the derivatives we get immediately
  \begin{equation}
  \label{Eq:3r}
  \partial_{X_r}\dogston=\frac{1}{X_r}\sum_{i=r}^ne_iX_\twodots{i}=\frac{1}{X_r}
  \left(\dogston-\sum_{i=0}^{r-1}e_iX_\twodots{i}\right)
  \end{equation}
  \begin{equation}
  \label{Eq:4r}
    \partial_{X_r}\dogstonwok=\frac{1}{X_r}\sum_{i=r}^{n-1}e_i^{(k)}X_\twodots{i}+
    \frac{1}{X_kX_r}\sum_{i=k}^{n-1}e_i^{(k)}X_\twodots{i}(X_{i+1}-X_k)
  \end{equation}
  \begin{equation}
  \label{Eq:5r}
    \hphantom{\partial_{X_r}\dogstonwok}=\frac{1}{X_r}\left(\dogstonwok-\sum_{i=0}^{r-1}e_i^{(k)}X_\twodots{i}\right)
  \end{equation}
  By plugging (\ref{Eq:3r}) and (\ref{Eq:5r}) into (\ref{Eq:0r}) we obtain
\[
X_r\Delta_r =\dogston\left( \sum_{i=0}^{r-1}e_i^{(k)}X_\twodots{i}\right)-
                            \dogstonwok\left( \sum_{i=0}^{r-1}e_iX_\twodots{i}\right)=
\]
and after simple cancelation, by invoking (\ref{Eq:2r}) we get
\[
\begin{array}{l}
=\left(\sum_{j=r}^{n}e_jX_\twodots{j}\right)\left(\sum_{i=0}^{r-1}e_i^{(k)}X_\twodots{i}\right)-\\
\hspace{5ex}
\left(\sum_{j=r}^{n-1}e_j^{(k)}X_\twodots{j}+
\frac{1}{X_k}\sum_{j=k}^{n-1}e_j^{(k)}X_\twodots{j}(X_{j+1}-X_k)\right)
\left(\sum_{i=0}^{r-1}e_iX_\twodots{i}\right)
\end{array}
\]
i.e.
{\small
  \[
    X_r\Delta_r=\!\!\!\!\!\!\sum_{0\leq i< r\leq j\leq n}\!\!\!\!\!\!(e_je_i^{(k)}-
    e_ie_j^{(k)})X_\twodots{i}X_\twodots{j}
    +\frac{1}{X_k}\sum_{0\leq i<r,k\leq j<n}\!\!\!\!\!\!\!\!e_ie_j^{(k)}X_\twodots{i}X_\twodots{j}(X_k-X_{j+1})
  \]
}
  If we use a simple identity $e_j=e_j^{(k)}+\xi_ke_{j-1}^{(k)}$, we can identify the quantity
  \[
  \begin{array}{l}
    e_je_i^{(k)}-e_ie_j^{(k)}=(e_j^{(k)}+\xi_ke_{j-1}^{(k)})e_i^{(k)}-(e_i^{(k)}+\xi_ke_{i-1}^{(k)})e_j^{(k)}=\\[2ex]
    \hphantom{e_je_i^{(k)}-e_ie_j^{(k)}}
    =\left|
    \begin{array}{cc}
    e_{j-1}^{(k)}&e_j^{(k)}\\
    e_{i-1}^{(k)}&e_i^{(k)}
    \end{array}
    \right|\xi_k
    =s^{(k)}_{2^i1^{j-i-1}}\xi_k
  \end{array}
  \]
  Thus in this case $(1\leq r<k)$ we obtain a formula
  \[
    \begin{array}{ll}
    \Delta_r= &\displaystyle \xi_k\!\!\!\!\!\!\!\!\!\!\sum_{0\leq i< r\leq j\leq n}\!\!\!\!\!\!
    s^{(k)}_{(j-1,i)^{(k)}}X_1^2\cdots X_i^2X_{i+1}\cdots \widehat{X}_k\cdots X_j+\\
    &\displaystyle\hphantom{xxxx}+\!\!\!\!\!\!\!\!\!\!\sum_{0\leq i<r,k\leq j<n}\!\!\!\!\!\!\!\!
    e_ie_j^{(k)}X_1^2\cdots X_i^2X_{i+1}\cdots \widehat{X}_r\cdots \widehat{X}_k\cdots X_j(X_k-X_{j+1})
    \end{array}
  \]
  (where $e_j^{(k)}=e_j^{(k)}=e_j(\xi_1,\ldots ,\widehat{\xi_k},\ldots,\xi_n)$)
  in terms of a Schur function (of arguments $\xi_1,\ldots ,\widehat{\xi_k},\ldots,\xi_n$)
  corresponding to a transpose $(2^i1^{j-i-1})$ of a partition $(j-i,i)$
  (cf.\ Jacobi--Trudi formula, I 3.5 in \cite{Jacobi}).\\[5mm]


  {\bf (ii)} For any $r$, $(1\leq )k<r\leq n$. In this case we use
  \[
\partial_{X_r}\dogstonwok=\frac{1}{X_kX_r}\sum_{j=r-1}^{n-1}e_j^{(k)}X_\twodots{j+1}
  \]
\[
\begin{array}{l@{\ =\ }l}
\dogstonwok
&\displaystyle
\sum_{i=0}^{k-1}e_i^{(k)}X_\twodots{i}+\frac{1}{X_k}\sum_{i=k}^{n-1}e_i^{(k)}X_\twodots{i+1}=\\
&\displaystyle
\frac{1}{X_k}\left(\sum_{i=0}^{k-1}X_\twodots{i}(X_k-X_{i+1})+\sum_{i=0}^{n-1}e_i^{(k)}X_\twodots{i}\right)
\end{array}
\]
By plugging this into (\ref{Eq:0r}) we get
\[
\begin{array}{ll}
X_kX_r\Delta_r
&\displaystyle =\left( \sum_{j=r}^{n}e_jX_\twodots{j}\right)
 \left( \sum_{i=0}^{k-1}e_i^{(k)}X_\twodots{i}(X_k-X_{i+1})+\sum_{i=0}^{n-1}e_i^{(k)}X_\twodots{i+1}\right)-\\
 &\displaystyle \hphantom{xxxx}-\left(\sum_{j=0}^{r-1}e_jX_\twodots{j}+\sum_{j=r}^{n}e_jX_\twodots{j}\right)
 \left(\sum_{i=r-1}^{n-1}e_i^{(k)}X_\twodots{i+1}\right)\\
&\hspace{-10ex}\displaystyle
 =\left( \sum_{i=0}^{r-2}e_i^{(k)}X_\twodots{i+1}\right)\left( \sum_{j=r}^{n}e_jX_\twodots{j}\right)-
\left( \sum_{i=0}^{r-1}e_iX_\twodots{i}\right)\left(\sum_{j=r-1}^{n-1}e_j^{(k)}X_\twodots{j+1} \right)+\\
 &\displaystyle \hphantom{xxxx}+\sum_{i=0}^{k-1}\sum_{j=r}^{n}e_i^{(k)}e_jX_\twodots{i}X_\twodots{j}(X_k-X_{i+1})\\
&\hspace{-10ex}\displaystyle
=\left( \sum_{i=1}^{r-1}e_{i-1}^{(k)}X_\twodots{i}\right)\left( \sum_{j=r}^{n}e_jX_\twodots{j}\right)-
\left( \sum_{i=0}^{r-1}e_iX_\twodots{i}\right)\left(\sum_{j=r}^{n}e_{j-1}^{(k)}X_\twodots{j} \right)+\\
 &\displaystyle \hphantom{xxxx}+\sum_{i=0}^{k-1}\sum_{j=r}^{n}e_i^{(k)}e_jX_\twodots{i}X_\twodots{j}(X_k-X_{i+1})\\
\end{array}
\]
By using a formula for elementary symmetric functions ($e_i=e_i^{(k)}+\xi_ke_{i-1}^{(k)}$)
we can write in terms of Schur functions (of arguments
$\xi_1,\ldots ,\xi_{k-1},\xi_{k+1},\ldots,\xi_n$), where ${\lambda '}$ is a conjugate
of $\lambda$.
\[
e_{i-1}^{(k)}e_j-e_ie_{j-1}^{(k)}=e_{i-1}^{(k)}e_j^{(k)}-e_i^{(k)}e_{j-1}^{(k)}=
-\left|
\begin{array}{cc}
e_{j-1}^{(k)} & e_j^{(k)}\\
e_{i-1}^{(k)} & e_i^{(k)}
\end{array}
\right|=-s^{(k)}_{2^i1^{j-i-1}}=-s^{(k)}_{(j-1,i)'}
\]
Thus we obtain a formula
\[
\begin{array}{ll}
\Delta_r=&\displaystyle -\left(
\sum_{0\leq i< r\leq j\leq n}\!\!\!\!\!\!
    s^{(k)}_{(j-1,i)'}X_1^2\cdots X_i^2X_{i+1}\cdots \widehat{X}_k\cdots \widehat{X}_r\cdots X_j+\right.\\
    &\displaystyle\hphantom{xxxx}\left. +\!\!\!\!\!\!\!\!\!\!\sum_{0\leq i<k,r\leq j<n}\!\!\!\!\!\!\!\!
    e_i^{(k)}e_jX_1^2\cdots X_i^2X_{i+1}\cdots \widehat{X}_k\cdots \widehat{X}_r\cdots X_j(X_{j+1}-X_k)
\right)
    \end{array}
\]
\begin{corollary}
\label{C:2.1}
Let $X_1\geq\cdots\geq X_n\geq 0$, $\xi_1,\ldots,\xi_n\geq 0$ be as before. Then
\begin{enumerate}[(i)]
  \item for any $r$, $1\leq r<k\ (\leq n)$ we have
  \[
\frac{\Psi^{1\ldots n}_{1\ldots n}}{
{\Psi^{1\ldots \widehat{k}\ldots n}_{1\ldots \widehat{k}\ldots n}}}
\geq
\frac{\Psi^{1\ldots\ r+1\ r+1\ \ldots n}_{1\ldots\ \hspace{1.2ex}r\hspace{2ex}r+1\ \ldots n}}{
{\Psi^{1\ldots\ r+1\ r+1\ \ldots\widehat{k}\ldots n}_{1\ldots\ \hspace{1.2ex}r\hspace{2ex}r+1\ \ldots \widehat{k}\ldots n}}}
  \]

  \item for any $r$, $(1\leq)\ k<r\ (\leq n)$ we have
  \[
\frac{\Psi^{1\ldots n}_{1\ldots n}}{
{\Psi^{1\ldots \widehat{k}\ldots n}_{1\ldots \widehat{k}\ldots n}}}
\geq
\frac{\Psi^{1\ldots\ r-1\ r-1\ \ldots n}_{1\ldots\ r-1\hspace{2ex}r\ \hspace{1.2ex}\ldots n}}
{{\Psi^{1\ldots\widehat{k}\ldots\ r-1\ r-1\ \ldots n}_{1\ldots\widehat{k}\ldots\ r-1\hspace{2ex}r\ \hspace{1.2ex}\ldots n}}}
  \]
\end{enumerate}
\end{corollary}
Now we illustrate how to use Corollary \ref{C:2.1} to prove our Conjecture \ref{C:1} for $n=2, 3, 4$ and $5$.\\
{\bf Case $n=2$}
\[
Q_2:=\frac{\dogsa{12}}{\dogsa{1}\dogsa{2}}\geq\frac{\dogs{22}{12}}{\dogs{2}{1}\dogsa{2}}=1\mbox{ (by $(i)$)}
\]
{\bf Case $n=3$}
\[
\begin{array}{rl}
\displaystyle
Q_3:=\frac{\dogsa{123}\dogsa{123}}{\dogsa{12}\dogsa{13}\dogsa{23}}&\displaystyle
\geq
\frac{\dogs{223}{123}\dogsa{123}}{\dogs{22}{12}\dogsa{13}\dogsa{23}}
\geq
\frac{\dogs{223}{123}\dogs{223}{123}}{\dogs{22}{12}\dogsa{13}\dogsa{23}}\mbox{ (by $(i)$)}\\[5mm]
&\displaystyle
\geq
\frac{\dogs{222}{123}\dogs{223}{123}}{\dogs{22}{12}\dogs{22}{13}\dogsa{23}}
\geq
\frac{\dogs{222}{123}\dogs{222}{123}}{\dogs{22}{12}\dogs{22}{13}\dogsa{23}}=1\mbox{ (by $(ii)$)}
\end{array}
\]
{\bf Case $n=4$}
\[
\begin{array}{rl}
\displaystyle
Q_4:=\frac{(\dogsa{1234})^3}{\dogsa{123}\dogsa{124}\dogsa{134}\dogsa{234}}
\geq\cdots\geq
\frac{\dogs{2244}{1234}(\dogs{2224}{1234})^2}{\dogs{224}{123}\dogs{224}{124}\dogs{224}{134}\dogs{224}{234}}
\ \ (\geq 1)
\end{array}
\]
This last inequality follows from the following symmetric function identity:
\[
\begin{array}{l}
\dogs{2244}{1234}(\dogs{2224}{1234})^2-\dogs{224}{123}\dogs{224}{124}\dogs{224}{134}\dogs{224}{234}=\\[2mm]
X_2^2X_4^4m_{2222}+2X_2^2X_4^3m_{2221}+X_2^2X_4^2m_{222}+3X_2^2X_4^2m_{2211}+X_2^2X_4m_{221}\\[2mm]
+4X_2^2X_4m_{2111}+X_2^2m_{211}+X_2(3X_2+2X_4)m_{1111}+X_2m_{111}
\end{array}
\]
where $m_\lambda=m_\lambda(\xi_1,\xi_2,\xi_3,\xi_4)$ are the monomial symmetric functions.\\
{\bf Case $n=5$}
\[
\begin{array}{rl}
\displaystyle
Q_5:=\frac{(\dogstona{5})^4}{\prod_{k=1}^5\dogstonwoka{5}}
\geq\cdots\geq
\frac{(\dogs{22244}{12345}\dogs{22444}{12345})^2}
{\dogs{2244}{1234}\dogs{2244}{1235}\dogs{2244}{1245}\dogs{2244}{1345}\dogs{2244}{2345}}
\ \ (\geq 1)
\end{array}
\]
The last inequality is equivalent to an explicit symmetric function identity with all
coefficients (w.r.t.\ monomial basis) positive.

Now we state our stronger conjecture.
\begin{conjecture}(for symmetric functions)
\begin{enumerate}[(a)]
\item For $n$ even we have
\[
\!\!\!\!\!\!\!\!\!\!\!\!
\dogs{2\ 2\ 4\ 4\ldots n\ n}{1\ 2\hspace{0.61ex}\ldots\hspace{0.61ex} n-1\ n}
\left(
\prod_{k=1}^{n/2}\dogs{2\ 2\ 4\ 4\ldots 2k\ 2k\ 2k\ldots n-2\ n-2\ n}
{1\ 2\ 3\ 4 \hspace{6.5ex}\ldots\hspace{7.2ex} n-1\ n}
\right)^2
\geq
\prod_{k=1}^n\dogs{2\ 2\ 4\ 4\ldots n-2\ n-2\ n}
{1\ 2 \hspace{1.31ex}\ldots\hspace{1.31ex} \widehat{k}\hspace{1.31ex}\ldots\hspace{1.31ex} n-1\ n}
\]

\item For $n$ odd we have
\[
\left(
\prod_{k=1}^{\lfloor n/2\rfloor}\dogs{2\ 2\ 4\ 4\ldots 2k\ 2k\ 2k\ldots n-1\ n-1}
{1\ 2\ 3\ 4\hspace{2.35ex}\ldots\hspace{6.8ex} n-1\hspace{2ex} n}
\right)^2
\geq
\prod_{k=1}^n\dogs{2\ 2\ 4\ 4\ldots n-1\ n-1}
{1\ 2 \hspace{1.57ex}\ldots\hspace{1.57ex} \widehat{k}\hspace{1.57ex}\ldots\hspace{1.57ex} n}
\]
\end{enumerate}
\end{conjecture}
Now we motivate another inequalities for symmetric functions which also refine the strongest
Atiyah--Sutcliffe conjecture for configurations of type (A). Let $n=3$.
We apply Corollary \ref{C:2.1} by using steps $(ii)$ only.
\[
Q_3:=\frac{\dogsa{123}\dogsa{123}}{\dogsa{12}\dogsa{13}\dogsa{23}}
\geq
\frac{\dogs{113}{123}\dogsa{123}}{\dogsa{12}\dogsa{13}\dogs{13}{23}}
\geq
\frac{\dogs{112}{123}\dogsa{123}}{\dogsa{12}\dogs{12}{13}\dogs{13}{23}}
\geq
\frac{\dogs{112}{123}\dogs{122}{123}}{\dogsa{12}\dogs{12}{13}\dogs{12}{23}}
\geq
1
\]
The last inequality is equivalent to nonnegativity of the expression
\[
\dogs{112}{123}\dogs{122}{123}-\dogsa{12}\dogs{12}{13}\dogs{12}{23}\ \
(=X_1(X_1-X_2)^2\xi_1\xi_2\xi_3\geq 0).
\]
Similarly, for $n=4$, the symmetric function inequality stronger than $Q_4\geq 1$
would be the following
\[
\dogs{1123}{1234}\dogs{1223}{1234}\dogs{1233}{1234}
\geq
\dogsa{123}\dogs{123}{124}\dogs{123}{134}\dogs{123}{234}
\]
Now we state a general conjecture for symmetric functions which imply the strongest
Atiyah--Sutcliffe conjecture for almost collinear type (A) configurations.
\begin{conjecture}
\label{C:2.2}
Let $X_1\geq\cdots\geq X_n\geq$, $\xi_1,\ldots\xi_n\geq 0$. Then the following inequality
for symmetric functions in $\xi_1,\ldots,\xi_n$
\[
\dogs{112\ldots n-1}{123\ldots n}
\dogs{1223\ldots n-1}{1234\ldots n}
\cdots
\dogs{12\ldots n-2\ n-1\ n-1}{12\ldots n-2\ n-1\ n}
\geq
\dogsa{1\ 2\ldots n-1}
\dogs{1\ 2\ldots n-1}{1\ 2\ldots n-2\ n}
\cdots
\dogs{1\ 2\ldots n-1}{2\ 3\ldots n-1}
\]
i.e.
\[
\prod_{k=1}^{n-1}
\dogs{1\ 2\ldots k\ \hphantom{+}k\hphantom{1}\ldots n}{1\ 2\ldots k\ k+1\ldots n}
\geq
\prod_{k=1}^{n}
\dogs{1\ 2\ \ldots\  n-1}{1\ 2\ldots \widehat{k}\ldots n}
\]
holds true.
\end{conjecture}

We have checked this Conjecture \ref{C:2.2} up to $n=5$ by using MAPLE and symmetric
function package of J.\ Stembridge. For $n$ bigger than five the computations are extremely
intensive and hopefully in the near future would be possible by using more powerful
computers.

Note that the right hand side of the Conjecture \ref{C:2.2} involves symmetric functions
of partial alphabets $\xi_1, \xi_2,\ldots,\xi_{k-1},\xi_{k+1},\ldots,\xi_n$. But the left hand side
doesn't have this "defect". Our objective now is to give explicit formula for the right hand
side in terms of the elementary symmetric functions of the full alphabet
$\xi_1, \xi_2,\ldots,\xi_n$. This we are going to achieve by using resultants as follows.
\begin{lemma}
\label{L:res}
For any $k$, ($1\leq k\leq n$), we have
\[
\dogs{1\ldots k\ldots n}{1\ldots \widehat{k}\ldots n}=
\sum_{j=0}^{n-1}a_j\xi_k^{n-1-j}
\]
where
\[
\begin{array}{l}
a_{n-1} = 1+X_1e_1+X_1X_2e_2+\ldots +X_1\cdots X_{n-1}e_{n-1},\\
a_{n-2} = -X_1-X_1X_2e_1-\ldots -X_1\cdots X_{n-1}e_{n-2},\\
\cdots\\
a_0 = (-1)^{n-1}X_1\cdots X_{n-1}\\[5mm]
\mbox{i.e.}\\
\displaystyle a_{n-1-j}=(-1)^j\sum_{i=j}^{n-1}X_1\cdots X_ie_{i-j}
\end{array}
\]
\end{lemma}
\proof[of Lemma \ref{L:res}]
By definition we have
\begin{equation}
\label{E:2.8}
\dogs{1\ldots n-1}{1\ldots\widehat{k}\ldots n}=\sum_{i=0}^{n-1}X_1\cdots X_ie_i^{(k)}
\end{equation}
where $e_i^{(k)}$ is the $i$--th elementary function of $\xi_1,\ldots ,\xi_{k-1},\xi_{k+1},\ldots,\xi_n$.
Now from the decomposition
\[
(1+\xi_kt)^{-1}\prod_{j=1}^{n}(1+\xi_jt)=\prod_{j\neq k}(1+\xi_jt)=
\sum_{i=0}^{n-1}e_i^{(k)}t^i
\]
we get
\[
e_i^{(k)}=e_i-e_{i-1}\xi_k+e_{i-2}\xi_k^2-\cdots +(-1)^{i}\xi_k^i
\]
By substituting this into equation (\ref{E:2.8}) the Lemma \ref{L:res} follows.
\qed\\

Then, by Lemma \ref{L:res}, the right hand side
\[
R_n=\prod_{k=1}^{n}\dogs{1\ 2\ \ldots\ k\ \ldots \  n-1}{1\ 2\ \ldots\  \widehat{k}\ \ldots\ \ \ n}
=
\prod_{k=1}^n\left(\sum_{j=0}^{n-1}a_j\xi_k^{n-1-j} \right)
\]
can be understood as a resultant $R_n=Res(f,g)$ of the following two polynomials
\[
\begin{array}{r@{\ =\ }l}
f(x) & \displaystyle \sum_{j=0}^{n-1}a_jx^{n-1-j}\\[2mm]
g(x) & \displaystyle \prod_{i=1}^{n}(x-\xi_i)=\sum_{j=0}^{n}(-1)^je_jx^{n-j}
\end{array}
\]
The Sylvester formula
\[
R_n=
\left|
\begin{array}{cccccccc}
  1 & -e_1 & e_2 & -e_3 & \ldots & (-1)^{n}e_n &   &  \\
    & 1 & -e_1 & e_2 & -e_3 & \ldots &   &  \\
    &   & \ddots &   &   &   &   &   \\
    &   &   & 1 &  -e_1 & \cdots  &   &   \\
  a_0 & a_1 & a_2 & \cdots & a_n &   &   &   \\
  & a_0 & a_1 & a_2 & \cdots & a_n &   &  \\
    &   & \ddots &   &   &   &   &   \\
    &   &   & a_0 & a_1 & a_2 & \cdots & a_n \\
\end{array}
\right|
\ \
\left(=:
\left|
\begin{array}{cc}
A & B\\
C & D
\end{array}
\right|
\right)
\]
can be simplified as
\[
=|A|\cdot |D-CA^{-1}B|=|D-CA^{-1}B|.
\]
The entries of the $n\times n$ matrix $\Delta:=D-CA^{-1}B$ are given by
\[
\delta_{ij}=\left\{
\begin{array}{ll}
\displaystyle \sum_{k=0}^{i-j}(-1)^{i+j}X_1\cdots X_{i-k}e_{j-k} & \mbox{ for } n-1\geq i\geq j\geq 0\\
\displaystyle \sum_{k=0}^{n}(-1)^{j-i}X_1\cdots X_{j-i+k}e_{j-i+k+1} & \mbox{ for } n-1\geq j> i\geq 0
\end{array}
\right.
\]
For example, for $n=3$
\[
\Delta_3=\left|
\begin{array}{ccc}
  1 & X_1e_2+X_1X_2e_3 & -X_1e_3 \\[2mm]
  -X_1 & 1+X_1e_1 & X_1X_2e_3 \\[2mm]
  X_1X_2 & -X_1-X_1X_2e_1 & 1+X_1e_1+X_1X_2e_2
\end{array}
\right|
\]
By elementary operations we get
\[
\Delta_3=
\left|
\begin{array}{ccc}
  1 & * & * \\[2mm]
  0 & \dogs{112}{123} & X_1(X_2-X_1)e_3 \\[2mm]
  0 & X_2-X_1 & \dogs{122}{123}
\end{array}
\right|
=
\left|
\begin{array}{cc}
\dogs{112}{123} & X_1(X_2-X_1)e_3 \\[2mm]
X_2-X_1 & \dogs{122}{123}
\end{array}
\right|
\]
Similarly, for $n=4$ we obtain
{\small
\[
\Delta_4=\left|
\begin{array}{ccc}
  \dogs{1123}{1234} & -X_1(X_1-X_2)e_3-X_1X_2(X_1-X_3)e_4 & X_1(X_1-X_2)e_4 \\[2mm]
  -(X_1-X_2) & \dogs{1223}{1234} & -X_1X_2(X_2-X_3)e_4 \\[2mm]
  X_1(X_2-X_3) & -(X_1-X_3)-X_1(X_2-X_3)e_1 & \dogs{1233}{1234}
\end{array}
\right|
\]
}
In general
\[
\Delta_n=\det(\delta'_{ij})_{1\leq i,j\leq n-1}
\]
where
\[
\delta'_{ij}=\left\{
\begin{array}{l@{\ ,\ }l}
\displaystyle \sum_{k=j+1}^{n-1}(-1)^{i+j+1}X_1\cdots X_{k-1}(X_k-X_i)e_{k+1} & \mbox{ for } i<j\\[2mm]
\displaystyle \dogs{1\ \ldots\ i\ i\ \ldots\ n-1}{1\ 2\ \ldots\ n} &  \mbox{ for } i=j\\[2mm]
\displaystyle \sum_{k=0}^{j-1}(-1)^{i+j}X_1\cdots X_{i-k-2}(X_{i-k-1}-X_i)e_{k} & \mbox{ for } i>j
\end{array}
\right.
\]
\begin{corollary}
The conjecture \ref{C:2.2} is equivalent to a Hadamard type
inequality for the (non Hermitian) matrix $\displaystyle
(\delta'_{ij})_{1\leq i,j\leq n-1}$, i.e.
\[
\det(\delta'_{ij})\leq \prod_{i=1}^{n-1}\delta'_{ii}.
\]
\end{corollary}

\end{document}